# Simulating time-harmonic acoustic wave effects induced by periodic holes/inclusions on surfaces

Wen Hu[1,2], Zhuojia Fu[1,*], Leevan Ling[2]

[1]Center for Numerical Simulation Software in Engineering and Sciences, College of Mechanics and Materials, Hohai University, Nanjing 211100, PR China

[2]Department of Mathematics, Hong Kong Baptist University, Kowloon Tong, Hong Kong

*Corresponding author: zhuojiafu@gmail.com (Zhuojia Fu)

**Abstract**

This paper introduces the first attempt to employ a localized meshless method to analyze time-harmonic acoustic wave propagation on curved surfaces with periodic holes/inclusions. In particular, the generalized finite difference method is used as a localized meshless technique to discretize the surface gradient and Laplace-Beltrami operators defined extrinsically in the governing equations. An absorbing boundary condition is introduced to reduce reflections from boundaries and accurately simulate wave propagation on unclosed surfaces with periodic inclusions. Several benchmark examples demonstrate the efficiency and accuracy of the proposed method in simulating acoustic wave propagation on surfaces with diverse geometries, including complex shapes and periodic holes or inclusions.

## 1. Introduction

In recent decades, there has been growing interest in studying partial differential equations (PDEs) on curved surfaces due to their applications in science and engineering. Examples include texture generation in computer graphics [1-5], water wave propagation on Earth's surface [6], and Turing spot pattern generation in biology [7-9], among others. Therefore, the development of numerical approaches for solving surface PDEs become a prominent research topic.

Since all surface differential operators are locally defined on the tangent space of the surface rather than in Euclidean space, traditional numerical approaches are unsuitable for directly solving surface PDEs. Appropriate reformulations and techniques are needed to deal with surface PDEs. They can be classified into three categories: intrinsic [10-13], embedding [14, 15] and extrinsic [16-18]. Intrinsic methods are suitable for solving surface partial differential equations (PDEs) because they define and compute differential operators purely in terms of quantities intrinsic to the surface itself. These methods work by defining variables in an atlas of the surface, which contains local parameterizations, a partition of unity, and a metric tensor that defines the inner product on the surface. Using these elements, intrinsic methods can locally parameterize parts of the surface and solve the PDEs on the parameterized domains. The local solutions can then be combined using the partition of unity to obtain a global solution on the entire surface. The metric tensor allows differential operators to be defined and PDEs formulated and solved intrinsically, without reference to an ambient Euclidean space. In comparison, embedding methods reformulate the surface PDEs so they can be solved in the ambient Euclidean space using standard numerical approaches. Extrinsic methods employ the analytical relationship between surface differential operators and those in Euclidean space to directly solve surface PDEs. While embedding and extrinsic methods can be more direct, intrinsic methods are appealing because they define and compute everything purely in terms of quantities intrinsic to the surface itself. However, parameterizing the surface and constructing the partition of unity can be complicated and computationally expensive for intrinsic methods.

By employing the aforementioned reformulation or technique, the finite element method (FEM) [19-23] is utilized to study diffusion dynamics, hydrodynamics on stationary and evolving surfaces, and bulk-surface coupling. However, generating high-quality meshes for complex geometries and remeshing evolving surfaces continue to be time-consuming and challenging in FEM. In contrast, meshless methods like the localized radial basis function collocation method (LRBFCM) [24-31], the generalized finite difference method (GFDM) [17, 32-38], and the RBF-based finite difference method (RBF-FD) [39-41] have been developed to solve PDEs on stationary and evolving surfaces. These methods only require node discretization, eliminating the need for remeshing and mitigating mesh distortion issues encountered in FEM.

Despite the utilization of numerical schemes to investigate various physical and mechanical phenomena on curved surfaces, acoustic wave propagation on curved surfaces has received limited research attention. This simulation is crucial for describing the propagation of sound waves on membranes and thin shell structures. This study presents the first attempt to simulate time-harmonic acoustic wave propagation on curved surfaces with periodic holes/inclusions using the GFDM to solve surface PDE in extrinsic form. This paper is organized as follows: Section 2 introduces the mathematical models for time-harmonic acoustic wave propagation on surfaces with periodic holes/inclusions. Section 3 presents the GFDM coupled with extrinsic for the mathematical models given in Section 2. Section 4 demonstrates the efficiency and accuracy of the proposed numerical scheme through numerical results from several benchmark examples. Finally, Section 5 provides the concluding remarks.

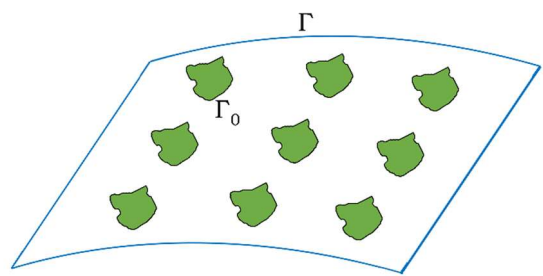

**Figure. 1.** Sketch of the surface with periodic inclusions (the blue line is the boundary of the surface).

## 2. Time-harmonic acoustic wave propagation model on the surfaces

Similar to phonon crystals [42-44], we study a smooth surface $S$ with periodic inclusions, as shown in Figure 1. The computational domain on $S$ consists of the matrix domain $D_1$ (white region) and the inclusion domain $D_2$ (green region), satisfying $D_1 \cup D_2 = S$, and $D_1 \cap D_2 = \Gamma_0$. The transient acoustic wave propagation on both the matrix and inclusions can be described by [45, 46],

$$\frac{\partial^2 \Phi^j}{\partial t^2} = c_j^2 \Delta_S \Phi^j + f_0, \bm{x} \in D_j, \tag{1}$$

where $\Phi^j = \Phi^j(x_1, x_2, x_3, t)(j=1,2)$ denotes the velocity potential, $c_j$ stands for the speed of acoustic wave propagating on surfaces, in which the superscript $j$ represents the physical quantities of either $D_1$ or $D_2$. And $f_0$ is the source term, $\Delta_S$ represents Laplace-Beltrami operator in tangent space, which is defined as,

$$\Delta_S \Phi^j := \nabla_S \Phi^j \cdot \nabla_S \Phi^j, \tag{2}$$

where the surface gradient $\nabla_S$ has the following form,

$$\nabla_S \Phi^j := (I - \vec{n}\vec{n}^T)\nabla \Phi^j, \tag{3}$$

and $\vec{n} = (n_1, n_2, n_3)^T$ represents the normal vector of the considered node $\bm{x} = (x_1, x_2, x_3) \in S$. The velocity potential and source term can be expressed as $\Phi^j(x_1, x_2, x_3, t) = Re\{u^j(x_1, x_2, x_3)e^{-i\omega t}\}$ and $f_0 = Re\{g_0(x_1, x_2, x_3)e^{-i\omega t}\}$, respectively, using the method of separation of variables. Subsequently, the governing equation (4) can be simplified to the following form,

$$\mu_j \Delta_S u^j + \rho_j \omega^2 u^j = g, \bm{x} \in D_j, \tag{4}$$

where $u^j(x_1, x_2, x_3)$ represents the displacement in $D_j$. $\mu_j$, $\rho_j$ and $\omega$ respectively represent shear modulus, mass density, and angular frequency, and $g = -\rho_j g_0$. On the interface boundaries $\Gamma_0$, the following continuity boundary conditions should be satisfied,

$$u^1_{\Gamma_0} = u^2_{\Gamma_0}, \mu_1 \partial u^1 / \partial \vec{n}_{\Gamma_0} = \mu_2 \partial u^2 / \partial \vec{n}_{\Gamma_0}. \tag{5}$$

The excitation of the incident acoustic wave is imposed on the boundary $\Gamma_I$ as,

$$u^1 = C, \boldsymbol{x} \in \Gamma_I, \tag{6}$$

where $C$ is an arbitrary constant. To reduce spurious reflections from the boundaries, an absorbing boundary condition (ABC) [45, 47] is applied on boundary $\Gamma_A$,

$$\frac{\partial u^1}{\partial \vec{\boldsymbol{n}}} - i\frac{\omega}{c}u^1 = 0, x \in \Gamma_A. \tag{7}$$

Like the phonon crystals, the surface $S$ is studied with semi-infinite periodic inclusions in the $x_p$ ($p=1,2,3$)-direction. the periodic boundaries $\Gamma_{c1}$ and $\Gamma_{c2}$ ($\boldsymbol{x}_{\Gamma_{p2}} = \boldsymbol{x}_{\Gamma_{p1}} + \boldsymbol{a}$) should satisfy Bloch theorem [48-50],

$$u^1_{\Gamma_{p2}} = e^{i\boldsymbol{k}\boldsymbol{a}} u^1_{\Gamma_{p1}}, \mu_1 \partial u^1 / \partial \vec{\boldsymbol{n}}_{\Gamma_{p2}} = e^{i\boldsymbol{k}\boldsymbol{a}} \mu_1 \partial u^1 / \partial \vec{\boldsymbol{n}}_{\Gamma_{p2}}, \tag{8}$$

where $\boldsymbol{k}$ is the Bloch wave vector, $\boldsymbol{a} = m_1\boldsymbol{a}_1 + m_2\boldsymbol{a}_2$, in which $\boldsymbol{a}_1$ and $\boldsymbol{a}_2$ express the fundamental translation vectors of the lattice, $m_1$ and $m_2$ are arbitrary integers. $\vec{\boldsymbol{n}}_{\Gamma_{pi}}, i=1,2$ represents the unit normal vector along the boundaries $\Gamma_{p1}$ or $\Gamma_{p2}$, and $\Gamma_P = \Gamma_{p1} \cup \Gamma_{p2}$. The remaining boundaries $\Gamma_N$ satisfy the zero Neumann boundary conditions,

$$\partial u^1 / \partial \vec{\boldsymbol{n}} = 0, \boldsymbol{x} \in \Gamma_N. \tag{9}$$

Hence, the surface boundaries can be given by $\Gamma = \Gamma_I \cup \Gamma_A \cup \Gamma_P \cup \Gamma_N$. By setting the governing equation (4), interface condition (5) and boundary conditions (6) - (9), one gets the following matrix equations,

$$AU = \begin{bmatrix} u^1_{\Gamma_I} & 0 \\ u^1_{\Gamma_{p2}} - e^{ika}u^1_{\Gamma_{p1}} & 0 \\ \vec{n}_{\Gamma_A} \cdot \nabla u^1 - i\dfrac{\omega}{c_1}u^1 & 0 \\ \mu_1 \partial u^1/\partial \vec{n}_{\Gamma_{p2}} - e^{ika}\mu_1 \partial u^1/\partial \vec{n}_{\Gamma_{p1}} & 0 \\ \partial u^1/\partial \vec{n}_{\Gamma_N} & 0 \\ u^1_{\Gamma_0} & -u^2_{\Gamma_0} \\ \mu_1 \partial u^1/\partial \vec{n}_{\Gamma_0} & -\mu_2 \partial u^2/\partial \vec{n}_{\Gamma_0} \\ (\mu_1 \Delta_S + \rho_1 \omega^2)u^1 & 0 \\ 0 & (\mu_2 \Delta_S + \rho_2 \omega^2)u^2 \end{bmatrix} \begin{bmatrix} C \\ 0 \end{bmatrix} = \begin{bmatrix} 0 \\ 0 \\ 0 \\ 0 \\ 0 \\ 0 \\ 0 \\ 0 \\ 0 \end{bmatrix} = b. \qquad (10)$$

## 3. GFDM for surface PDEs

This section presents an extrinsic meshless method for solving the surface Helmholtz equations. The extrinsic form scheme represents the Laplace-Beltrami operator and surface gradient operator in tangent spaces as the Laplace operator and gradient operator in Euclidean space.

By referring to the surface operators definitions in (2) and (3) and a theorem in Ref. [15], the Laplace-Beltrami operator can be expressed as,

$$\Delta_S u^j := \Delta u^j - H_S \partial_{\vec{n}} u^j - \partial_{\vec{n}}^{(2)} u^j, \qquad (11)$$

where $\partial_{\vec{n}} u^j := \boldsymbol{n}^T \nabla u^j$, $\partial_{\vec{n}}^{(2)} u := \vec{\boldsymbol{n}} J(\nabla u)\vec{\boldsymbol{n}}$ and $H_S(\boldsymbol{x}) := tr(J(\vec{\boldsymbol{n}})(I - \vec{\boldsymbol{n}}\vec{\boldsymbol{n}}^T))$, and $J$ is the Jacobian operator in Euclidean space. By substituting definition (11) into the governing equation (4), the following equations can be derived,

$$\mu_j \left[ -H_S(x_1, x_2, x_3)\left(n_1 \frac{\partial u^j}{\partial x_1} + n_2 \frac{\partial u^j}{\partial x_2} + n_3 \frac{\partial u^j}{\partial x_3}\right) + (1-n_1^2)\frac{\partial^2 u^j}{\partial x_1^2} + (1-n_2^2)\frac{\partial^2 u^j}{\partial x_2^2} \right. \\ \left. + (1-n_3^2)\frac{\partial^2 u^j}{\partial x_3^2} - 2\left(n_1 n_2 \frac{\partial^2 u^j}{\partial x_1 \partial x_2} + n_1 n_3 \frac{\partial^2 u^j}{\partial x_1 \partial x_3} + n_2 n_3 \frac{\partial^2 u^j}{\partial x_2 \partial x_3}\right) \right] + \rho_j \omega^2 u^j = g. \qquad (12)$$

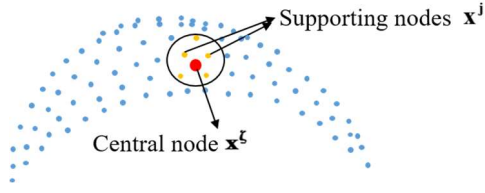

**Figure. 2.** Schematic diagram of node distribution on the surface.

The meshless GFDM is employed to discretize PDE (12). In Figure 2, each node in the computational domain forms a subdomain by selecting $m$ nodes closest to the central node $\mathbf{x}^\zeta$. The GFDM is implemented by performing a second-order Taylor series expansion at $\mathbf{x}^\zeta$,

$$u_j = u_\zeta + h_j \frac{\partial u}{\partial x_1}\bigg|_{\mathbf{x}^\zeta} + k_j \frac{\partial u}{\partial x_2}\bigg|_{\mathbf{x}^\zeta} + l_j \frac{\partial u}{\partial x_3}\bigg|_{\mathbf{x}^\zeta} + \frac{h_j^2}{2} \frac{\partial^2 u}{\partial x_1^2}\bigg|_{\mathbf{x}^\zeta} + \frac{k_j^2}{2} \frac{\partial^2 u}{\partial x_2^2}\bigg|_{\mathbf{x}^\zeta} + \frac{l_j^2}{2} \frac{\partial^2 u}{\partial x_3^2}\bigg|_{\mathbf{x}^\zeta}$$
$$+ h_j k_j \frac{\partial^2 u}{\partial x_1 \partial x_2}\bigg|_{\mathbf{x}^\zeta} + h_j l_j \frac{\partial^2 u}{\partial x_1 \partial x_3}\bigg|_{\mathbf{x}^\zeta} + k_j l_j \frac{\partial^2 u}{\partial x_2 \partial x_3}\bigg|_{\mathbf{x}^\zeta} + \cdots, \quad j = 1, 2, \ldots m, \quad (13)$$

where $u_\zeta, u_j$ express the solution at the central node $\mathbf{x}^\zeta$ and the supporting nodes $\{\mathbf{x}^j\}_{j=1}^m$, and $h_j$, $k_j$, $l_j$ are the coordinate components of distance between $\mathbf{x}^\zeta$ and $\mathbf{x}^j$,

$$h_j = x_1^j - x_1^\zeta, k_j = x_2^j - x_2^\zeta, l_j = x_3^j - x_3^\zeta. \quad (14)$$

Next, the residual function $R(u)$ [51] in each subdomain can be defined as follows,

$$R(u) = \sum_{j=1}^m \left[ \left( u_\zeta - u_j + h_j \frac{\partial u}{\partial x_1}\bigg|_{\mathbf{x}^\zeta} + k_j \frac{\partial u}{\partial x_2}\bigg|_{\mathbf{x}^\zeta} + l_j \frac{\partial u}{\partial x_3}\bigg|_{\mathbf{x}^\zeta} + \frac{h_j^2}{2} \frac{\partial^2 u}{\partial x_1^2}\bigg|_{\mathbf{x}^\zeta} + \frac{k_j^2}{2} \frac{\partial^2 u}{\partial x_2^2}\bigg|_{\mathbf{x}^\zeta} + \frac{l_j^2}{2} \frac{\partial^2 u}{\partial x_3^2}\bigg|_{\mathbf{x}^\zeta} \right.\right.$$
$$\left.\left. + h_j k_j \frac{\partial^2 u}{\partial x_1 \partial x_2}\bigg|_{\mathbf{x}^\zeta} + h_j l_j \frac{\partial^2 u}{\partial x_1 \partial x_3}\bigg|_{\mathbf{x}^\zeta} + k_j l_j \frac{\partial^2 u}{\partial x_2 \partial x_3}\bigg|_{\mathbf{x}^\zeta} \right) w_j \right]^2, \quad (15)$$

in which the weighting function $w_j$ can be chosen as the following fourth-order spline function [51],

$$w_j = \begin{cases} 0, & d_j > d_{max} \\ 1 - 6\left(\frac{d_j}{d_{max}}\right)^2 + 8\left(\frac{d_j}{d_{max}}\right)^3 - 3\left(\frac{d_j}{d_{max}}\right)^4, & d_j \leq d_{max} \end{cases}, \quad (16)$$

where $\{d_j\}_{j=1}^m$ represent the distances from the central node $\mathbf{x}^\zeta$ to the supporting nodes $\mathbf{x}^j$, and $d_{max}$ stands for the maximum distance among all the distances $\{d_j\}_{j=1}^m$.

By minimizing the residual function with respect to partial derivatives

($\frac{\partial u}{\partial x_1}\big|_{\mathbf{x}^\iota}$, $\frac{\partial u}{\partial x_2}\big|_{\mathbf{x}^\iota}$, ..., $\frac{\partial^2 u}{\partial x_2 \partial x_3}\big|_{\mathbf{x}^\iota}$), the following discretization scheme for the partial differential terms in GFDM can be derived,

$$\frac{\partial u}{\partial x_i}\bigg|_{\mathbf{x}^\iota} = w_0^{x_i} u_\zeta + \sum_{j=1}^{m} w_j^{x_i} u_j, \quad \frac{\partial^2 u}{\partial x_i \partial x_k}\bigg|_{\mathbf{x}^\iota} = w_0^{x_i x_k} u_\zeta + \sum_{j=1}^{m} w_j^{x_i x_k} u_j, \quad i,k=1,2,3, \qquad (17)$$

After substituting discretization scheme (17) into PDE (12), the governing equation can be expressed in matrix form as follows,

$$A_{N \times N} U_{N \times 1} = B, \qquad (18)$$

where $N$ represents the total number of computational points in the surface domain. $A_{N \times N}$ is the sparse coefficient matrix obtained using the discretization scheme (17), with each row having $m$ non-zero elements.

## 4. Numerical results and discussions

In this section, we numerically investigate examples of acoustic wave propagation on surfaces. To evaluate the accuracy of the results, we define the global relative error as follows,

$$\text{Global Error} = \left[ \sum_{i=1}^{N_T} (u_{numerical}^i - u_{exact}^i)^2 \right]^{1/2} \bigg/ \left[ \sum_{i=1}^{N_T} (u_{exact}^i)^2 \right]^{1/2}, \qquad (19)$$

where $N_T$ stands for the test nodes number on the considered surface, $u_{numerical}^i$ and $u_{exact}^i$ represent the numerical and exact solutions at the test node $\mathbf{x}^i$, respectively. Next, the acoustic wave propagation behaviors on the surfaces with complex geometries and periodic inclusions are numerically investigated.

Unless stated otherwise, $N_T = N$, the supporting node numbers is set to $m=40$, the constant boundary condition (6) is chosen as $C = 10^{-5}$ for the system balance. And the materials used for the matrix and inclusions are epoxy and aurum (Au) in this study. The acoustic velocities and densities of epoxy and Au are known as $c_0 = 1161 m/s$, $\rho_0 = 1180 kg/m^3$, $c_1 = 1239 m/s$, $\rho_1 = 19,500 kg/m^3$.

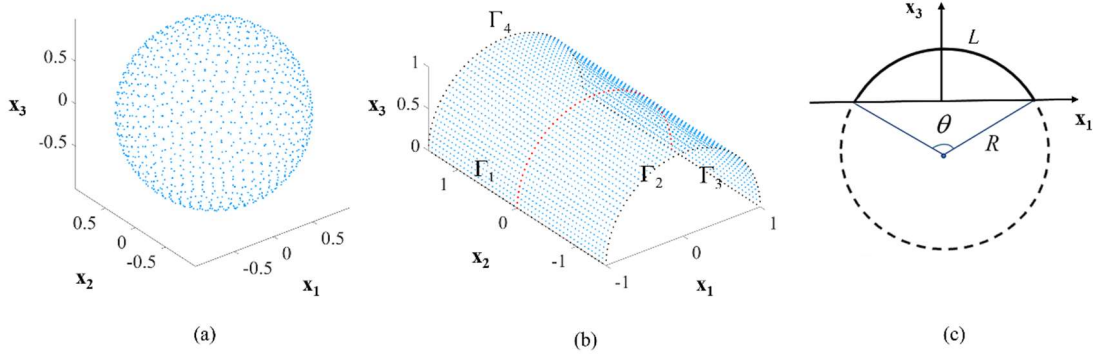

**Figure. 3.** Sketch of node distributions on (a) closed spherical surface, (b) portion of the cylinder surface; (c) side view of the portion of the cylinder surface.

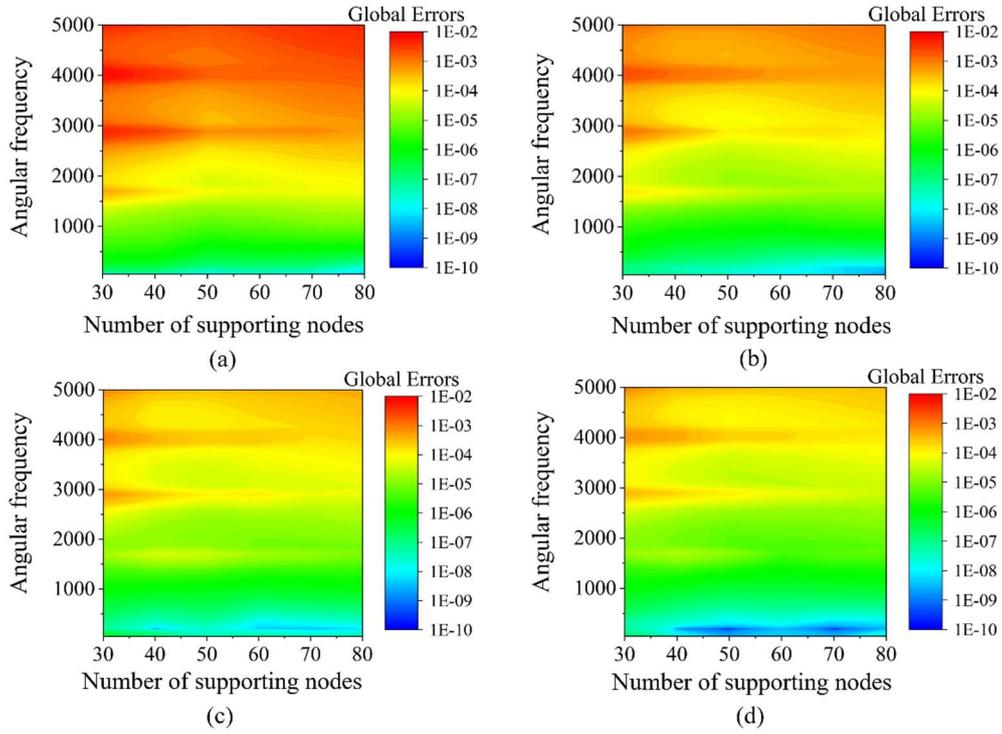

**Figure. 4.** Global error variations with respect to angular frequencies $\omega$ and supporting node numbers $m$ by using the proposed method with (a) $N$=2500, (b) $N$=4900, (c) $N$=8100, (d) $N$=10000 on the closed spherical surface.

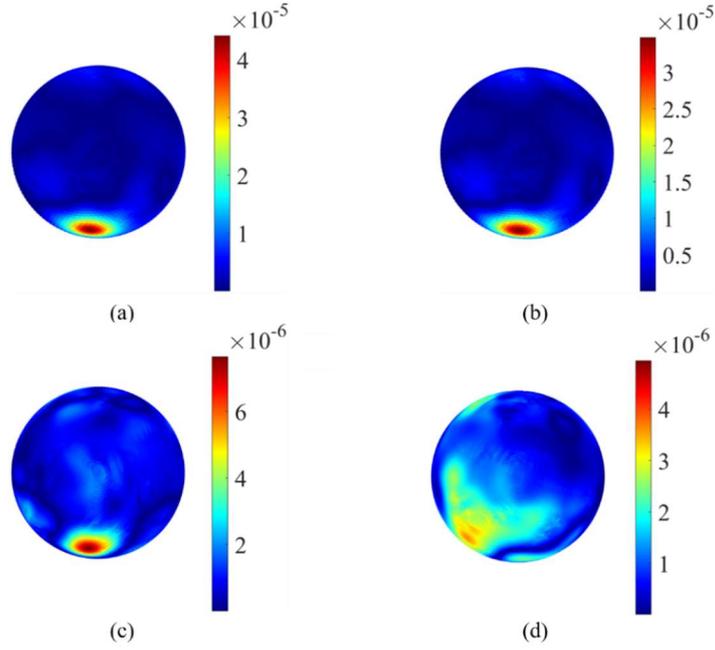

**Figure. 5.** Relative error distributions on the closed spherical surface obtained by using the proposed method with various total node numbers: (a) $N = 2500$; (b) $N = 4900$; (c) $N = 8100$; (d) $N = 10000$.

4.1 Convergence and numerical efficiency analysis

This section validated the accuracy, convergence, and efficiency of the proposed method for acoustic wave propagation. We examine two types of surfaces: a closed spherical surface and a simple unclosed surface, both made of epoxy, i.e., $D_1 = S, D_2 = \varnothing$.

Case 1.1 Closed spherical surface

For the closed spherical surface ($S_{a1} = x_1^2 + x_2^2 + x_3^2 - 1$) shown in Figure 3(a), consider the governing equation (4) with the boundary $\Gamma = \varnothing$. The source term g in the governing equation (4) is expressed as,

$$g = 10^{-5} \frac{\omega^2}{c_1^2} \left(2(x_1 x_2 + x_1 x_3 + x_2 x_3) - 1\right) \cos\left(\frac{\omega(x_1 + x_2 + x_3)}{c_1}\right)$$
$$+ 2 \times 10^{-5} \frac{\omega}{c_1} (x_1 + x_2 + x_3) \sin\left(\frac{\omega(x_1 + x_2 + x_3)}{c_1}\right), (x_1, x_2, x_3) \in S_{a1}, \quad (20)$$

which can be easily derived from the following exact solution,

$$u^1 = 10^{-5} \cos\left(\frac{\omega(x_1 + x_2 + x_3)}{c_1}\right), (x_1, x_2, x_3) \in S_{a1}. \quad (21)$$

In this case, uniformly distributed nodes on the spherical surface are generated using the minimum energy scheme [17]. The proposed GFDM is applied to solve case 1.1 on the closed spherical surface with different node discretizations (N=2500, 4900, 8100, 10000).

Figure 4 illustrates global error variations for GFDM results with angular frequencies $\omega \in [50, 5000]$ and supporting node numbers $m \in [30, 80]$. The global relative errors are below $6.9 \times 10^{-3}$. Figure 5 displays relative error distributions with different node discretizations on the closed spherical surface (N=2500, 4900, 8100, 10000), using $m$=40 and $\omega$=1000. The results show two main trends. First, the numerical accuracy improves as the total node number or supporting node number increases, as seen in Figures 4 and 5. Secondly, the numerical errors decrease with decreasing angular frequency in the range $[50, 5000]$, as seen in Figures 4 and 5.

Case 1.2 Unclosed surface - a portion of the cylinder surface

Next, we consider governing equation (4) with zero source term ($g = 0$) and an angular frequency of $\omega$=10000, subjected to the boundary conditions (6) and (9) on the unclosed surface (a portion of a cylinder surface) defined as,

$$S_{a2} = \left\{(x_1, x_2, x_3) \middle| x_1^2 + x_3^2 = R^2, |x_1| \leq R \sin\left(\frac{\theta}{2}\right), |x_2| \leq \lambda, 0 \leq x_3 \leq R - R \cos\left(\frac{\theta}{2}\right), \theta \in (0, \pi]\right\},$$
(22)

where $\theta$ and $R$ represent the central angle corresponding to the arc length and the radius of the cylindrical cross-section, respectively. The arc length $L$ is given by $L = \theta R$, as illustrated in Figure 3(c). We set the parameters as $\lambda = 1.5$, $R$=1 and $\theta = \pi$, with the related nodes being uniformly distributed on $S_{a2}$, as depicted in Figure 3 (b). In this case, $\Gamma_I = \Gamma_1$, $\Gamma_A = \Gamma_P = \varnothing$, and $\Gamma_N = \Gamma_2 \cup \Gamma_3 \cup \Gamma_4$.

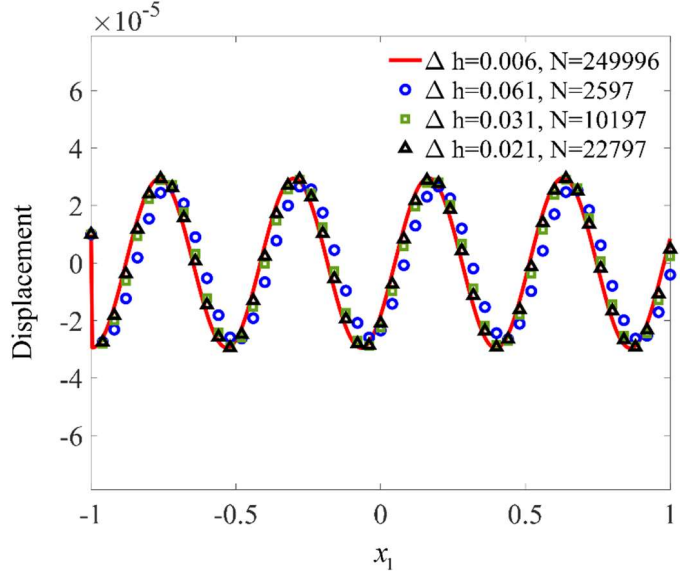

**Figure. 6.** Displacement variation on the test nodes ($N_T$=51) uniformly distributed on the cylindrical curve $\{(x_1, x_2, x_3) \in S_{a2} | -1 \leq x_1 \leq 1, x_2 = 0\}$.

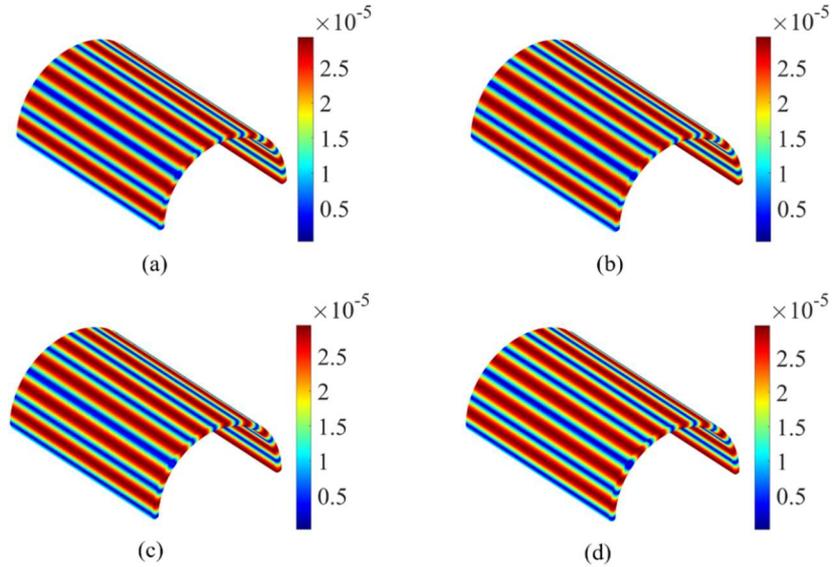

**Figure. 7.** Displacement distribution on the unclosed surface $S_{a2}$ obtained by using the proposed GFDM with different node densities: (a) $\Delta h = 0.061 (N = 2597)$; (b) $\Delta h = 0.031 (N = 10197)$; (c) $\Delta h = 0.021 (N = 22797)$; (d) Reference results: $\Delta h = 0.006 (N = 249996)$.

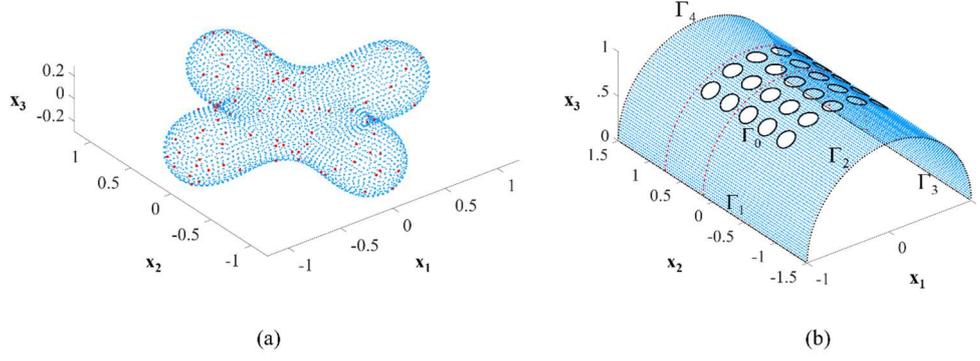

**Figure. 8.** Node distribution on the surfaces (a) $S_{b1}$; (b) $S_{b2}$.

Since there is no analytical solution available for this case, we utilize the proposed GFDM results with a dense node distribution of $\Delta h = 0.006$ ($N = 249996$) as reference solutions. Figure 6 presents the GFDM results on the test nodes (red dots in Figure 3(b)) uniformly distributed along the cylindrical curve $\{(x_1, x_2, x_3) \in S_{a2} | -1 \leq x_1 \leq 1, x_2 = 0\}$, using different node discretizations ($\Delta h = 0.061 (N = 2597)$; $\Delta h = 0.031 (N = 10197)$; $\Delta h = 0.021 (N = 22797)$). Figure 6 demonstrates the convergence of the GFDM results to the reference results as the node number increases. Figure 7 showcases the numerical displacement distributions on the unclosed surface $S_{a2}$. It is noteworthy that even with coarse nodes $\Delta h = 0.061 (N = 2597)$, the GFDM yields acceptable results.

4.2 Surfaces with complicated geometries

This section investigates acoustic wave propagation on closed and unclosed surfaces with complex geometries to showcase the efficiency of the proposed method. The surface geometries include a closed surface with a complex shape (constant distance product surface) and a portion of a cylindrical surface with periodic holes, representing an unclosed surface with complicated features.

Case 2.1 Closed constant distance product surface

In this case, we consider governing equation (4) with angular frequency of

$\omega = 10000$ on the following closed constant distance product surface,

$$S_{b1} = \sqrt{(x_1-1)^2 + x_2^2 + x_3^2}\sqrt{(x_1+1)^2 + x_2^2 + x_3^2}\sqrt{x_1^2 + (x_2-1)^2 + x_3^2}$$
$$\sqrt{x_1^2 + (x_2+1)^2 + x_3^2} - 1.1. \tag{23}$$

It has no boundary, i.e., $\Gamma = \varnothing$. The material of the entire surface $S_{b1}$ is specified as epoxy, namely, $D_1 = S, D_2 = \varnothing$, and the source term $g$ in governing equation (4) can be expressed as follows,

$$g = 10^{-5}\left(\left(\frac{\omega^2}{c_1^2} + n_1^2 - 1\right)\cos(x_1) + \left(\frac{\omega^2}{c_1^2} + n_2^2 - 1\right)\cos(x_2) + \left(\frac{\omega^2}{c_1^2} + n_3^2 - 1\right)\cos(x_3)\right.$$
$$\left. + H_s(x_1,x_2,x_3)(n_1 \sin x_1 + n_2 \sin x_2 + n_3 \sin x_3)\right), (x_1,x_2,x_3) \in S_{b1}, \tag{24}$$

which can be readily derived from the following exact solution,

$$u = 10^{-5}(\cos(x_1) + \cos(x_2) + \cos(x_3)), (x_1,x_2,x_3) \in S_{b1}. \tag{25}$$

In this example, the discretized nodes are placed on $S_{b1}$ (N=3996), as shown in Figure 8(a). Figure 9 presents the displacement variation on the test nodes (red dots in Figure 8(a), $N_T$=86) distributed on $S_{b1}$. Figure 10 displays the displacement distribution and the relative error distribution on the constant distance product surface $S_{b1}$. Figures 9 and 10 demonstrate good agreement between the numerical and exact displacements on $S_{b1}$, with a maximum relative error below $4\times10^{-5}$. This confirms the effectiveness of the proposed localized extrinsic meshless method for acoustic wave propagation on closed surfaces with complex geometries.

Case 2.2 Unclosed surface - a portion of the cylinder surface with periodic holes

We consider governing equation (4) with zero source term ($g = 0$) and angular frequency $\omega$=10000, subject to the boundary conditions (6)-(7) on the unclosed surface. A portion of a cylinder surface with 5×5 periodic circular holes as shown in Figure 8 (b), which can be defined as,

$$S_{b2} = \left\{(x_1,x_2,x_3) \in S_{a2} \,\middle|\, X_1^2 + (X_2 + 0.9 - 0.3i)^2 \geq 0.01, [X_1,X_2,X_3]^T = J_j[x_1,x_2,x_3]^T,\right.$$
$$\left. (i, j = 1,2,3,4,5)\right\}, \tag{26}$$

where $J_j(j = 1,2,3,4,5)$ represent the coordinate transformation matrix between

$[X_1, X_2, X_3]^T$ and $[x_1, x_2, x_3]^T$, and has the following form,

$$J_j = \begin{bmatrix} \cos\alpha_j & 0 & -\sin\alpha_j \\ 0 & 1 & 0 \\ \sin\alpha_j & 0 & \cos\alpha_j \end{bmatrix}, (\alpha_j = -0.9 + 0.3j, j = 1,2,3,4,5). \quad (27)$$

In this case, we set the parameters as $\lambda = 1.5$, $R = 1$ and $\theta = \pi$, and the nodes are uniformly distributed on $S_{b2}$, as shown in Figure 8 (b). In this case, $\Gamma_I = \Gamma_1$, $\Gamma_A = \Gamma_2 \cup \Gamma_3 \cup \Gamma_4$, $\Gamma_P = \Gamma_N = \varnothing$. Not completely the same as inclusions, the interface condition (5) between the matrix and holes can be rewritten as follows,

$$\mu_1 \partial u^1 / \partial \vec{n}_{\Gamma_0} = 0. \quad (28)$$

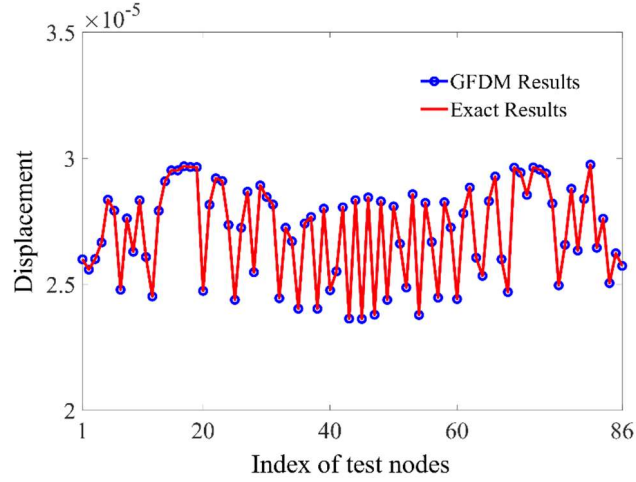

**Figure. 9.** Displacement variation on the test nodes (red dots in Fig. 8 (a), $N_T$=86) distributed on the constant distance product surface.

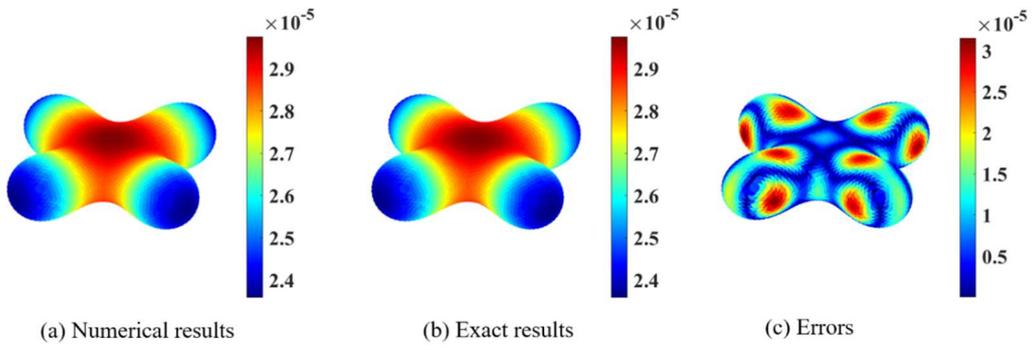

(a) Numerical results  (b) Exact results  (c) Errors

**Figure. 10.** Displacement distributions of (a) numerical solutions; (b) exact solutions; and relative error distributions on $S_{b1}$.

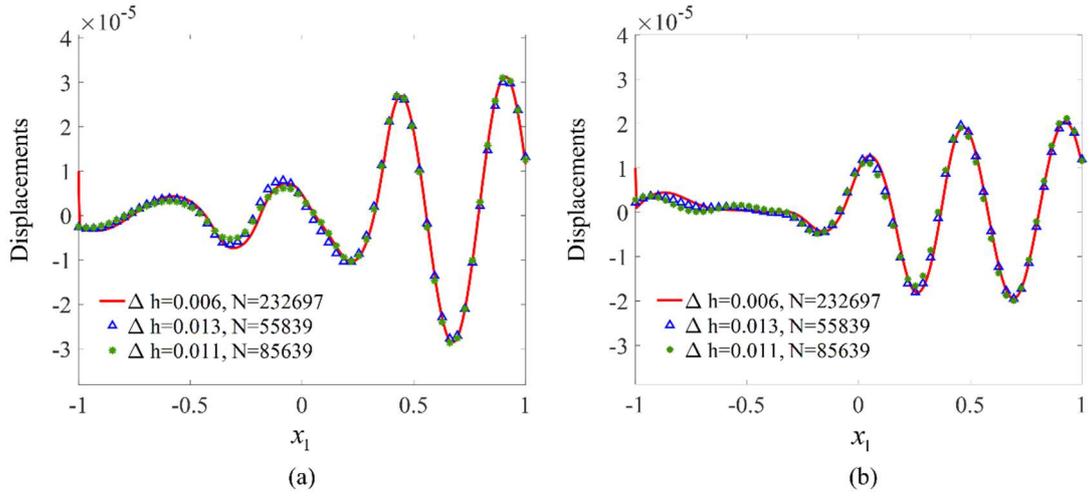

**Figure. 11.** Displacement variation on the test nodes (red dots in Fig. 8 (b)) distributed on the cylindrical curves, (a) $\{(x_1,x_2,x_3)\in S_{b2}|-1\leq x_1 \leq 1, x_2=0.15\}$, (b) $\{(x_1,x_2,x_3)\in S_{b2}|-1\leq x_1 \leq 1, x_2=0.75\}$.

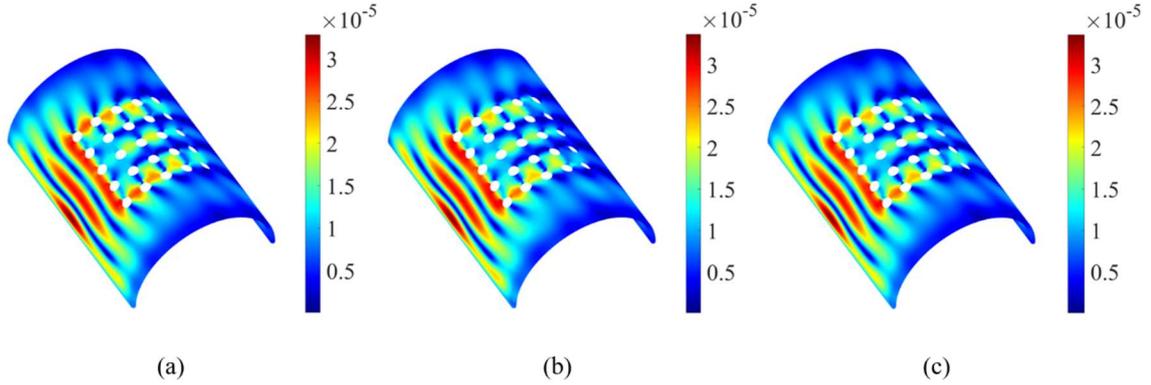

**Figure. 12.** Displacement distribution on the unclosed surface with 5×5 periodic circular holes $S_{b2}$ obtained by using the proposed GFDM with different node densities: (a) $\Delta h = 0.013(N=55839)$; (b) $\Delta h = 0.011(N=85639)$; (c) Reference results: $\Delta h = 0.006(N=232697)$.

Since there is no exact solution available for this case, we use the proposed GFDM solutions with a very dense node density of $\Delta h = 0.006(N=232697)$ as the reference solutions. Figure 11 depicts the GFDM results at the test nodes (red dots in Fig. 8(b)) distributed along the cylindrical curves $\{(x_1,x_2,x_3)\in S_{b2}|-1\leq x_1 \leq 1, x_2=0.15\}$ and $\{(x_1,x_2,x_3)\in S_{b2}|-1\leq x_1 \leq 1, x_2 = 0.75\}$ by using different node discretizations ( $\Delta h = 0.013(N=55839); \Delta h = 0.011(N=85639)$ ). Figure 11 demonstrate that the GFDM results converge as the node number increases. Figures 12 (a)-(c) present the numerical displacement distributions on $S_{b2}$ by using different node discretizations, and get the similar numerical results, indicating that the GFDM performs well on this surface with periodic holes.

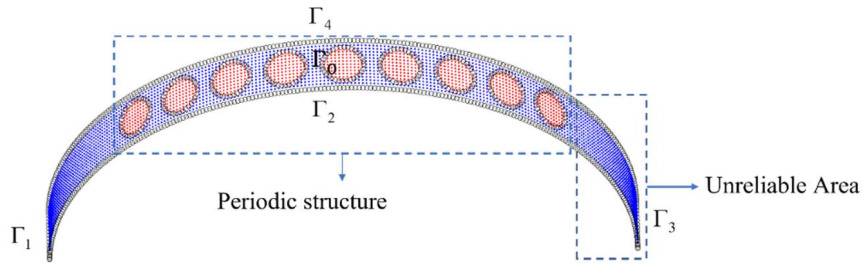

**Figure. 13.** Node distribution on the surface in Section 4.3.

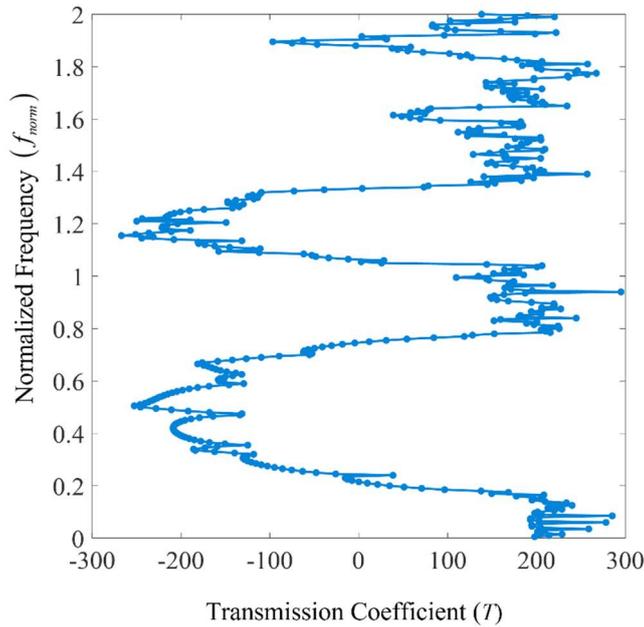

**Figure. 14.** Transmission spectra of aurum/epoxy materials on the unclosed surface $S_c$.

4.3 Surface with periodic inclusions

In this case study, we investigate an unclosed surface with semi-infinite periodic inclusions in the y-direction. To simplify this problem, only one inclusion with $S_c = S_{c1} \cup S_{c2}$ domain consists of an unclosed surface (epoxy, blue dots) with 9 circular inclusions (Au, red dots), as shown in Figure 13.

In this case, $D_1 = S_{c1}$, $D_2 = S_{c2}$, $\Gamma_I = \Gamma_1$, $\Gamma_A = \Gamma_3$, $\Gamma_P = \Gamma_2 \cup \Gamma_4$, and $\Gamma_N = \varnothing$. The matrix domain $S_{c1}$ and the circular inclusion domain $S_{c2}$ can be defined as follows,

$$S_{c1} = \{(x_1, x_2, x_3) | x_1^2 + x_3^2 = R^2, |x_1| \leq R, |x_2| \leq \lambda, 0 \leq x_3 \leq R, X_1^2 + X_2^2 \geq r_i^2,$$
$$[X_1, X_2, X_3]^T = J_j [x_1, x_2, x_3]^T, j = 1, 2, ..., 9\}, \quad (32)$$

$$S_{c2} = \{(x_1, x_2, x_3) | x_1^2 + x_3^2 = R^2, |x_1| \leq R, |x_2| \leq \lambda, 0 \leq x_3 \leq R, X_1^2 + X_2^2 \leq r_i^2,$$
$$[X_1, X_2, X_3]^T = J_j [x_1, x_2, x_3]^T, j = 1, 2, ..., 9\}, \quad (33)$$

in which $S_c = S_{c1} \cup S_{c2}$, $\Gamma_0 = S_{c1} \cap S_{c2}$, $\lambda = 0.5$, $R = \dfrac{16}{\pi}$ and $\theta = \pi$, $r_i$ stands for the radius of each circular inclusion, the coordinate transformation matrix $J_j$ can be given by using Eq. (27) with $\alpha_j = \dfrac{\pi}{16}(j-5), (j = 1, 2, ..., 9)$. Governed by equation (4), the related interface/boundary conditions (5)-(9) are imposed on the boundaries $\Gamma_0 - \Gamma_4$. Due to the Bloch theorem [48-50], as well as the periodicity and symmetry of the structure on $S_c$, the periodic boundary condition (8) imposed on $\Gamma_2$ and $\Gamma_4$ can be rewritten as follows,

$$u^1_{\Gamma_{p2}} = u^1_{\Gamma_{p1}}, \mu_1 \partial u^1 / \partial \vec{n}_{\Gamma_{p2}} = \mu_1 \partial u^1 / \partial \vec{n}_{\Gamma_{p2}}. \quad (29)$$

In this paper, the filling fraction and the normalized frequency are defined as $F_f = \pi r_i^2$, and $f_{norm} = \dfrac{\omega}{2\pi c_0}$, respectively. The transmission coefficient $T = 20 \times \log(u^0(\bm{x}_{\Gamma_3})/2 \times 10^{-5} \times u^0(\bm{x}_{\Gamma_1}))$ is used to describe the transmission spectra and exhibit the wave propagation characteristics on the surface. Figure 14 illustrates $T$ varying with the normalized frequency $f_{norm}$ in the range $[0.05, 2]$ (corresponding to the angular frequency range of $\omega \in [38, 14589]$) when $F_f = 0.4$ is maintained. Figure 14 illustrates that the transmission coefficient $T$ can be significantly small with normalized frequency $f_{norm}$ in a range of $[0.24, 0.74] \cup [1.06, 1.33]$. This indicates that acoustic waves with $f_{norm}$ among this specific range are unable to propagate through the periodic structures on the unclosed surface. Figure 15 presents the displacement distributions at two specific normalized frequencies: $f_{norm} = 0.5,\ 0.8$, to validate the abnormal acoustic propagation behaviors observed in Figure 14. As shown in Figure 15, the acoustic wave with $f_{norm} = 0.5$ is unable to propagate through the surface with periodic structures, while the acoustic wave with $f_{norm} = 0.8$ can pass through it, confirming the observations from Figure 14.

Similar to phonon crystals, the concept of a bandgap is introduced to represent the frequency range where the transmission coefficient $T < -10$. Figure 16 illustrates the effect of the filling fraction on the bandgap of the unclosed surface with periodic inclusions. In Figure 16, the first bandgap of the surface presents insensitivity to the filling fraction $F_f \leq 0.6$, while the second bandgap is highly influenced by the filling fraction within the range $F_f \in [0.2, 0.75]$.

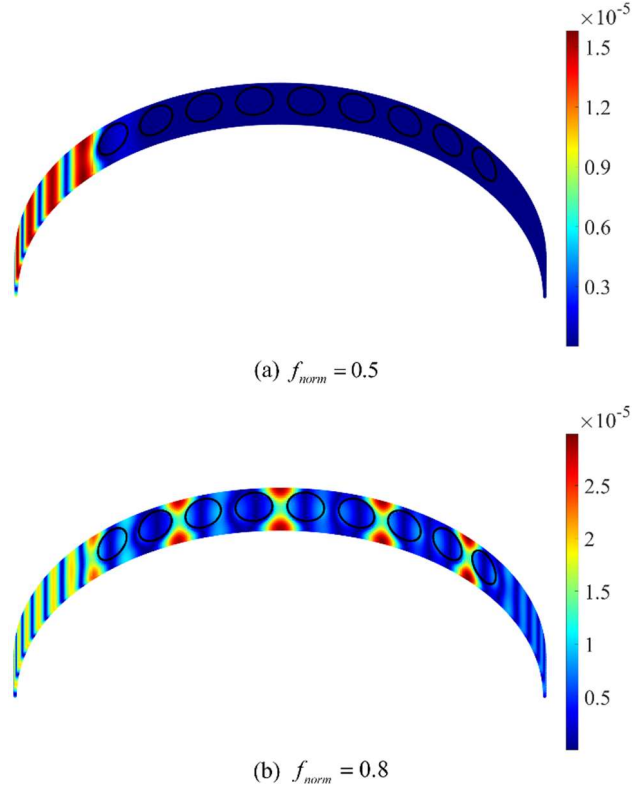

(a) $f_{norm} = 0.5$

(b) $f_{norm} = 0.8$

**Figure. 15.** Displacement amplitude distributions on the present unclosed surface $S_c$ generated by the acoustic wave with normalized frequencies.

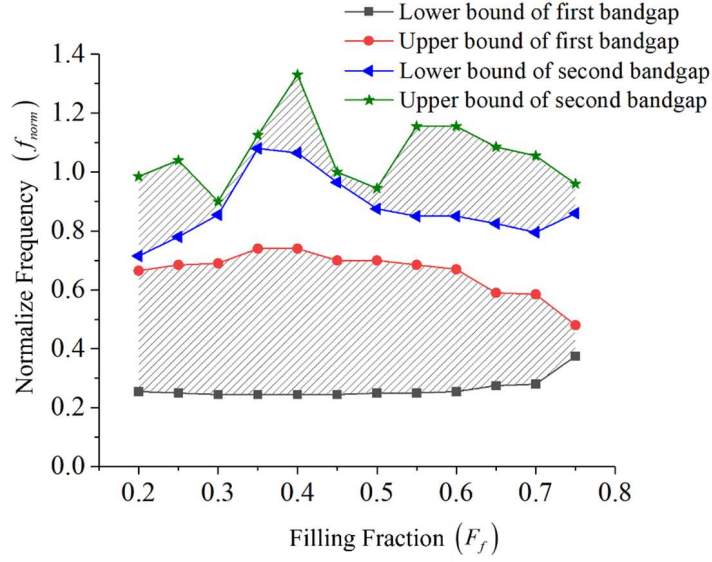

**Figure. 16.** Bandgap variation with respect to varied filling fractions $F_f \in [0.2, 0.75]$ on the present unclosed surface $S_c$

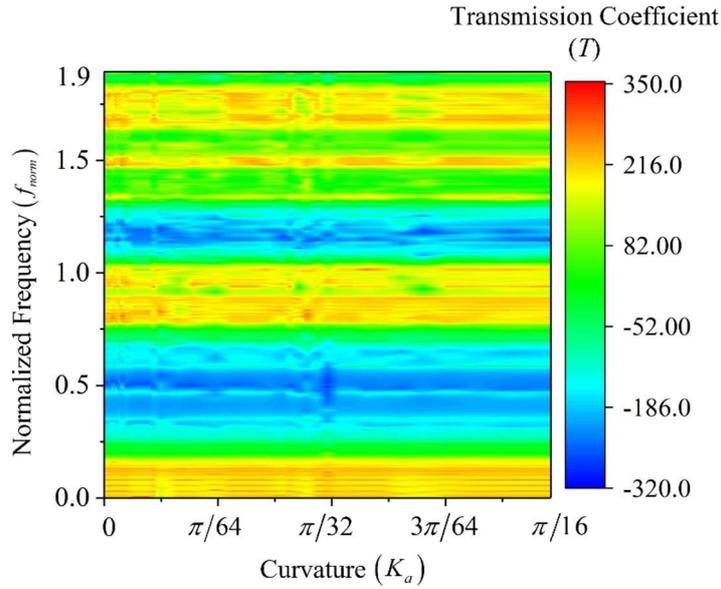

**Figure. 17.** The variation of the transmission coefficients $T$ with respect to varied curvatures $K_a \in [0, \pi/16]$ on the present unclosed surface $S_c$.

Then, we numerically investigate the effect of the surface curvature $K_a$ on $S_c$, where $K_a$ is defined as $1/R$ ($R$ being the radius of curvature). In the numerical implementation, the arc length of Sc is fixed at $L = 16$. Figure 17 illustrates the variation of the transmission coefficients $T$ on the unclosed surface $S_c$ with respect to different curvatures $K_a$. It can be observed that the curvature variation has a minimal effect on the bandgap region.

**Conclusions**

This paper presents an innovative approach for simulating time-harmonic acoustic wave propagation on curved surfaces. The methodology introduces an extrinsic GFDM coupled with an absorbing boundary condition to eliminate reflection effects. Numerical results demonstrate the efficiency and accuracy of the proposed method for simulating the acoustic wave propagation behavior on various closed surfaces with complex geometries, as well as unclosed surfaces with periodic holes or inclusions. As the number of total nodes and supporting nodes increases, the GFDM results converge towards the exact/reference results. Higher frequencies necessitate a larger number of nodes for achieving acceptable accuracy. Periodic holes/inclusions on curved surfaces can significantly impact the propagation of acoustic wave at specific frequencies. Finally, the inclusion filling fraction influences the bandwidth of the second bandgap, whereas curvature variation has a minimal impact on the bandgap region.

Additionally, this study focuses on wave propagation on stationary closed and unclosed surfaces with complex geometries and periodic holes or inclusions, as the initial step of our research. The proposed extrinsic GFDM for wave propagation on moving closed/unclosed surfaces, with periodic structures, is currently being intensively studied.

**Acknowledgments**

The work described in this paper was supported by the National Science Funds of China (Grant No. 12122205), the Six Talent Peaks Project in Jiangsu Province of China (Grant No. 2019-KTHY-009) and by the Hong Kong Research Grant Council GRF (Grant No. 12301520, 12301021, 12300922).

**References**

[1] J. Dorsey, P. Hanrahan, Digital materials and virtual weathering, Sci Am, 282 (2000) 64-71.
[2] Z. Yang, H. Liu, B. Zhang, W. Wang, W. Xu, Metasurfaces design for tuning of flexural wave and SH wave, Applied Physics A, 128 (2022) 695.
[3] A. Witkin, M. Kass, Reaction-diffusion textures, ACM SIGGRAPH Computer Graphics, 25 (1991) 299-308.

[4] H. Kim, A. Yun, S. Yoon, C. Lee, J. Park, J. Kim, Pattern formation in reaction–diffusion systems on evolving surfaces, Computers & Mathematics with Applications, 80 (2020) 2019-2028.

[5] C.M. Elliott, B. Stinner, C. Venkataraman, Modelling cell motility and chemotaxis with evolving surface finite elements, J R Soc Interface, 9 (2012) 3027-3044.

[6] T.G. Myers, J.P.F. Charpin, A mathematical model for atmospheric ice accretion and water flow on a cold surface, International Journal of Heat and Mass Transfer, 47 (2004) 5483-5500.

[7] M.A. Chaplain, M. Ganesh, I.G. Graham, Spatio-temporal pattern formation on spherical surfaces: numerical simulation and application to solid tumour growth, J Math Biol, 42 (2001) 387-423.

[8] I.F. Sbalzarini, A. Hayer, A. Helenius, P. Koumoutsakos, Simulations of (an)isotropic diffusion on curved biological surfaces, Biophys J, 90 (2006) 878-885.

[9] C.M. Elliott, B. Stinner, Modeling and computation of two phase geometric biomembranes using surface finite elements, J Comput Phys, 229 (2010) 6585-6612.

[10] M.S. Floater, K. Hormann, Surface Parameterization: a Tutorial and Survey, in: Advances in Multiresolution for Geometric Modelling, 2005, pp. 157-186.

[11] N. Trask, P. Kuberry, Compatible meshfree discretization of surface PDEs, Computational Particle Mechanics, 7 (2019) 271-277.

[12] A. Torres-Sánchez, D. Santos-Oliván, M. Arroyo, Approximation of tensor fields on surfaces of arbitrary topology based on local Monge parametrizations, J Comput Phys, 405 (2020) 109168.

[13] R. Lai, J. Liang, H.-K. Zhao, A local mesh method for solving PDEs on point clouds, Inverse Problems & Imaging, 7 (2013) 737-755.

[14] S.J. Ruuth, B. Merriman, A simple embedding method for solving partial differential equations on surfaces, J Comput Phys, 227 (2008) 1943-1961.

[15] K.C. Cheung, L. Ling, A Kernel-Based Embedding Method and Convergence Analysis for Surfaces PDEs, SIAM Journal on Scientific Computing, 40 (2018) A266-A287.

[16] M. Chen, L. Ling, Extrinsic Meshless Collocation Methods for PDEs on Manifolds, SIAM Journal on Numerical Analysis, 58 (2020) 988-1007.

[17] Z.C. Tang, Z.J. Fu, M. Chen, L. Ling, A localized extrinsic collocation method for Turing pattern formations on surfaces, Applied Mathematics Letters, 122 (2021) 107534.

[18] Z. Tang, Z. Fu, M. Chen, J. Huang, An efficient collocation method for long-time simulation of heat and mass transport on evolving surfaces, J Comput Phys, 463 (2022) 111310.

[19] W. Axmann, P. Kuchment, An Efficient Finite Element Method for Computing Spectra of Photonic and Acoustic Band-Gap Materials, J Comput Phys, 150 (1999) 468-481.

[20] A.S. Phani, J. Woodhouse, N.A. Fleck, Wave propagation in two-dimensional periodic lattices, J Acoust Soc Am, 119 (2006) 1995-2005.

[21] J. Grande, Eulerian Finite Element Methods for Parabolic Equations on Moving Surfaces, SIAM Journal on Scientific Computing, 36 (2014) B248-B271.

[22] M.A. Olshanskii, A. Reusken, X. Xu, A stabilized finite element method for advection-diffusion equations on surfaces, IMA Journal of Numerical Analysis, 34 (2013) 732-758.
[23] G. Dziuk, C.M. Elliott, Finite element methods for surface PDEs, Acta Numerica, 22 (2013) 289-396.
[24] H. Zheng, Z. Yang, C. Zhang, M. Tyrer, A local radial basis function collocation method for band structure computation of phononic crystals with scatterers of arbitrary geometry, Applied Mathematical Modelling, 60 (2018) 447-459.
[25] C.-M. Fan, C.-H. Yang, W.-S. Lai, Numerical solutions of two-dimensional flow fields by using the localized method of approximate particular solutions, Eng Anal Bound Elem, 57 (2015) 47-57.
[26] C.M. Fan, C.S. Chien, H.F. Chan, C.L. Chiu, The local RBF collocation method for solving the double-diffusive natural convection in fluid-saturated porous media, International Journal of Heat and Mass Transfer, 57 (2013) 500-503.
[27] B. Mavrič, D.N.M.P. Bozidar Sarler, B. Šarler, Local radial basis function collocation method for linear thermoelasticity in two dimensions, International Journal of Numerical Methods for Heat & Fluid Flow, 25 (2015) 1488-1510.
[28] I. Siraj ul, R. Vertnik, B. Šarler, Local radial basis function collocation method along with explicit time stepping for hyperbolic partial differential equations, Applied Numerical Mathematics, 67 (2013) 136-151.
[29] A. Karageorghis, C.S. Chen, Y.-S. Smyrlis, A matrix decomposition RBF algorithm: Approximation of functions and their derivatives, Applied Numerical Mathematics, 57 (2007) 304-319.
[30] C.S. Chen, A. Karageorghis, F. Dou, A novel RBF collocation method using fictitious centres, Applied Mathematics Letters, 101 (2020) 106069.
[31] M.A. Jankowska, A. Karageorghis, C.S. Chen, Improved Kansa RBF method for the solution of nonlinear boundary value problems, Eng Anal Bound Elem, 87 (2018) 173-183.
[32] Z.C. Tang, Z.J. Fu, H.G. Sun, X.T. Liu, An Efficient Localized Collocation Solver for Anomalous Diffusion on Surfaces, Fract Calc Appl Anal, 24 (2021) 865-894.
[33] J.J. Benito, A. García, L. Gavete, M. Negreanu, F. Ureña, A.M. Vargas, Solving a reaction–diffusion system with chemotaxis and non-local terms using Generalized Finite Difference Method. Study of the convergence, J Comput Appl Math, 389 (2021) 113325.
[34] Y. Gu, H.G. Sun, A meshless method for solving three-dimensional time fractional diffusion equation with variable-order derivatives, Applied Mathematical Modelling, 78 (2020) 539-549.
[35] Z.-J. Fu, W.-H. Chu, M. Yang, P.-W. Li, C.-M. Fan, Estimation of Tumor Characteristics in a Skin Tissue by a Meshless Collocation Solver, International Journal of Computational Methods, 18 (2020) 2041009.
[36] P.W. Li, Space–time generalized finite difference nonlinear model for solving unsteady Burgers' equations, Applied Mathematics Letters, 114 (2021) 106896.
[37] Z.-J. Fu, Z.-C. Tang, H.-T. Zhao, P.-W. Li, T. Rabczuk, Numerical solutions of

the coupled unsteady nonlinear convection-diffusion equations based on generalized finite difference method, The European Physical Journal Plus, 134 (2019) 272.

[38] P.W. Li, C.M. Fan, J.K. Grabski, A meshless generalized finite difference method for solving shallow water equations with the flux limiter technique, Eng Anal Bound Elem, 131 (2021) 159-173.

[39] A. Petras, L. Ling, C. Piret, S.J. Ruuth, A least-squares implicit RBF-FD closest point method and applications to PDEs on moving surfaces, J Comput Phys, 381 (2019) 146-161.

[40] A. Petras, L. Ling, S.J. Ruuth, An RBF-FD closest point method for solving PDEs on surfaces, J Comput Phys, 370 (2018) 43-57.

[41] B. Fornberg, E. Lehto, C. Powell, Stable calculation of Gaussian-based RBF-FD stencils, Computers & Mathematics with Applications, 65 (2013) 627-637.

[42] H.-W. Dong, S.-D. Zhao, X.-B. Miao, C. Shen, X. Zhang, Z. Zhao, C. Zhang, Y.-S. Wang, L. Cheng, Customized broadband pentamode metamaterials by topology optimization, Journal of the Mechanics and Physics of Solids, 152 (2021) 104407.

[43] H.-W. Dong, S.-D. Zhao, Y.-S. Wang, L. Cheng, C. Zhang, Robust 2D/3D multi-polar acoustic metamaterials with broadband double negativity, Journal of the Mechanics and Physics of Solids, 137 (2020) 103889.

[44] S.-D. Zhao, H.-W. Dong, X.-B. Miao, Y.-S. Wang, C. Zhang, Broadband Coding Metasurfaces with 2-bit Manipulations, Physical Review Applied, 17 (2022) 034019.

[45] Z.-J. Fu, A.-L. Li, C. Zhang, C.-M. Fan, X.-Y. Zhuang, A localized meshless collocation method for bandgap calculation of anti-plane waves in 2D solid phononic crystals, Eng Anal Bound Elem, 119 (2020) 162-182.

[46] C. Wang, L. Huang, Time-domain simulation of acoustic wave propagation and interaction with flexible structures using Chebyshev collocation method, J Sound Vib, 331 (2012) 4343-4358.

[47] Z.-J. Fu, L.-F. Li, D.-S. Yin, L.-L. Yuan, A Localized Collocation Solver Based on T-Complete Functions for Anti-Plane Transverse Elastic Wave Propagation Analysis in 2D Phononic Crystals, Mathematical and Computational Applications, 26 (2020) 2-22.

[48] J. Gazalet, S. Dupont, J.C. Kastelik, Q. Rolland, B. Djafari-Rouhani, A tutorial survey on waves propagating in periodic media: Electronic, photonic and phononic crystals. Perception of the Bloch theorem in both real and Fourier domains, Wave Motion, 50 (2013) 619-654.

[49] R.D.L. Kronig, W.G. Penney, Quantum mechanics of electrons in crystal lattices, Proceedings of the Royal Society of London. Series A, Containing Papers of a Mathematical and Physical Character, 130 (1997) 499-513.

[50] H. Watanabe, A Proof of the Bloch Theorem for Lattice Models, Journal of Statistical Physics, 177 (2019) 717-726.

[51] J.J. Benito, F. Urena, L. Gavete, Solving parabolic and hyperbolic equations by the generalized finite difference method, J Comput Appl Math, 209 (2007) 208-233.